\documentclass[12pt]{article}
\usepackage{mathrsfs}
\usepackage{}
\usepackage{amsfonts}
\usepackage{amssymb}
\usepackage{amsmath}
\usepackage{amsthm}
\usepackage{latexsym}
 \newcommand{\op} {\overline{\partial}}
\begin{document}
\numberwithin{equation}{section}
\title{Covering stability of Bergman kernels on K\"ahler hyperbolic manifolds}
\author{Xu Wang\footnote{Address: Department of Mathematics, Tongji University, Shanghai
200092, P. R. China. E-mail: 1113xuwang@mail.tongji.edu.cn}}
\date{}
\maketitle
\begin{abstract}
This paper is a sequel to \cite{Xu}. In this paper, an estimation of the Bergman Kernel of K\"ahler hyperbolic manifold is given by the $L^2$ estimate and the Bochner formula. As an application, an effective criterion of the very ampleness of the canonical line bundle of K\"ahler hyperbolic manifold is given, which is a generalization of Yeung's result.
\end{abstract}

\section{Introduction}

The notion of K\"ahler hyperbolic is due to Gromov \cite{Gromov91}. A \emph{non-compact} K\"ahler manifold $(X, \omega, \lambda)$ is called \emph{non-compact K\"ahler hyperbolic} if $\omega$ is the exterior differential of a $C^1$ bounded 1-form $\eta$, i.e., $\omega=d\eta$ and $|\eta|^2\leq\lambda$ on $X$ for some positive constant $\lambda$.
A compact K\"ahler manifold $(X, \omega, \lambda)$ is called \emph{K\"ahler hyperbolic} if the lift $\widetilde{\omega}$ of $\omega$ to the universal covering
$$
\widetilde{X}\rightarrow X=\widetilde{X}/\Gamma
$$
is the exterior differential of a $C^1$ bounded 1-form $\eta$, i.e., $\widetilde{\omega}=d\eta$ and $|\eta|^2_{\widetilde{\omega}}\leq \lambda$ on $\widetilde{X}$ (we say that $X$ is d-bounded by $\lambda$). What's more, if $\widetilde{X}$ is CH (Cartan-Hadamard) manifold (see Appendix II), we call $(X, \omega, \lambda)$ \emph{CH K\"ahler hyperbolic manifold}. Let $h^{n,0}=p_1$ be the first plurigenera, i.e., the dimension of the Bergman space of $X$. The $L_2$-Hodge number $\widetilde{h^{n,0}}=\widetilde{p_1}$ is defined as the integration of the Bergman kernel form of $\widetilde{X}$ on a Dirichlet fundamental domain of $X$ in $\widetilde{X}$.

Since $(X, \omega, \lambda)$ is compact, its Ricci curvature is bounded, \emph{throughout this paper, we shall assume that its Ricci curvature is bounded below by $-1$ unless specified mentioned}. Denote by $|X|$ the volume of $X$. Let
$\tau=\min_{x\in \widetilde{X}}\{ \tau(x) \}$, where $\tau(x)$ is the quasi-injectivity radius (see Appendix III). We shall prove that

\medskip

\textbf{Theorem 1.1.} {\it Let $(X, \omega, \lambda)$ be an $n$-dimensional CH K\"ahler hyperbolic manifold, if $n \geq 2$ and $\tau \geq 2\sqrt{2n}$, then}
\begin{equation}
\frac{1}{|X|}| \ p_1-\widetilde{p_1}|\leq 16(\frac{4}{\pi})^n\frac{\sqrt{n}\lambda}{\tau^2}.
\end{equation}

\medskip

For the compact ball quotients, by the same method as in the proof of Theorem 1.1, one could get a better estimation (see \cite{Xu} ). On the other hand, by a similar argument as in the proof of Theorem 1.1, we get the following generalization of Theorem 1.1.

\medskip

\textbf{Theorem 1.2.} {\it Let $(X, \omega, \lambda)$ be an $n$ dimensional K\"ahler hyperbolic manifold, if $\tau\geq 2$,}
\begin{equation}
\frac{1}{|X|}| \ p_1-\widetilde{p_1}|\leq 2^{203(2n+\sqrt{2n})}\frac{\lambda}{|B|\tau^2},
\end{equation}
{\it where $|B|$ stands for the minimal volume of the unit ball in $X$.}

\medskip

If $\widetilde X$ is a bounded homogeneous domain $\Omega$ in $\mathbb C^n$. By the result of Kai-Ohsawa \cite{Kai_Ohsawa07} (see also Vinberg-Gindikin-Pjatecki\u{\i}-\v{S}apiro \cite{V-G-P63}), one may choose suitable globally coordinate $z$ of $\Omega$ such that its Bergman kernel $K (z,\bar w)dz\otimes d\bar w$
satisfies
\begin{equation}
\widetilde{\omega} = i\partial\op \log\widetilde K(z,z), \ \ |\op \log \widetilde K(z,z)| \equiv KO_{\Omega},
\end{equation}
where $KO_{\Omega}$ is a positive constant only depends on the complex structure of $\Omega$. If
$$
\Omega=\{(u,v)\in\mathbb C^p \times \mathbb C^q ~|~ v+\bar v-F(u,u)\in V\}
$$
is a Siegel domain of second kind defined by $V$ and $F$, where $V$ is a convex cone in $\mathbb R^q$ containing no entire straight lines and $F$ is $V$-Hermitian. Ishi \cite{Ishi} proved that $KO_{\Omega}=\sqrt{p+2q}$.

\medskip

Since $\Omega$ is homogeneous, its Bergman kernel function $S_{\Omega}$ is a constant that only depends on the complex structure of $\Omega$. Thus its Ricci curvature satisfies
\begin{equation}
Ric(\widetilde{\omega}) = -\widetilde{\omega}.
\end{equation}
It is well known that the Bergman metric $\widetilde{\omega}$ on every bounded symmetric domain has non-positive sectional curvature. On the other hand, according to D'atri-Miatello \cite{DM83}, a bounded homogeneous domain with negative sectional curvature with respect to $\widetilde{\omega}$ must be symmetric. Thus, according to the above two Theorems, $p_1$ will be non-vanishing for sufficiently large $\tau$.

\medskip

\textbf{Theorem 1.3.} {\it The Bergman space of the compact quotient of a bounded symmetric domain $\Omega$ in $\mathbb C^n$ is nontrivial provided that}
\begin{equation}
\tau >  \max\{2^{n+1}\frac{KO_{\Omega}}{\sqrt{S_{\Omega}}} , \ 2\sqrt{2n}\}.
\end{equation}

\medskip

Since $\frac{1}{|X|}| \ p_1-\widetilde{p_1}|$ does't depend on the un-ramified normal covering of $X$, thus it could be used to investigate the covering behavior of the Bergman kernels. If a K\"ahler hyperbolic manifold $X$ has a tower of coverings $\{X_j\}$, by the result of Chen-Fu \cite{ChenFu} (see also \cite{Donnelly96}, \cite{Kazhdan70}, \cite{Ohsawa09}, \cite{Rhodes93}), the pull back of the Bergman kernel function on $X_j$ converges uniformly to that of $\widetilde{X}$ as $j\to \infty$. Thus it is natural to find an effective estimate of the Bergman kernel function. We shall prove an effective Ramadanov's Theorem (see \cite{Ramadanov67}) in section 3. Thus, we get an effective estimates of the Bergman kernels on K\"ahler hyperbolic manifolds, which is also given in section 3.

Inspired by Yeung's results (see \cite{Yeung00a}, \cite{Yeung00b}, \cite{Yeung94}), we shall give an effective ampleness criterion of the canonical line bundle in section 4 by using Calabi's diastasis function (see \cite{Calabi53}) and the Bergman metric instead of the heat kernel.

\section{The geometry of K\"ahler hyperbolic manifolds}

The following Theorem on the lower bound of the spectrum of the Laplace operator is due to Donnelly-Fefferman \cite{DonnellyFefferman83} and Gromov \cite{Gromov91} (see also Ohsawa \cite{Ohsawa09}). See Appendix-I for the basic notions (such as $\widetilde{\op}$ in (7.1)) of Hodge theory.

\medskip

\textbf{Theorem 2.1.} {\it Let $(X, \omega, \lambda)$ be a non-compact complete K\"ahler hyperbolic manifold of dimension $n$. If  Then every $u\in Dom \ \widetilde{\op}\cap Dom \ (\widetilde{\op})^*$ with degree $p+q\neq n$ satisfies the inequality}
\begin{equation}
||\widetilde{\op}u||^2 + ||(\widetilde{\op})^*u||^2 \geq \frac{|n-p-q|}{8\lambda}||u||^2.
\end{equation}
{\it For $u\in Dom \  \widetilde{\op}\cap Dom \ (\widetilde{\op})^*$ with degree $p+q = n$, which are orthogonal to the $L^2$ harmonic space,}
\begin{equation}
||\widetilde{\op}u||^2 + ||(\widetilde{\op})^*u||^2 \geq \frac{1}{8\lambda}||u||^2.
\end{equation}
{\it What's more, if $\omega$ has a global $C^2$ real potential $\psi$ such that}
\begin{equation}
\omega=i\partial\op\psi\geq\frac{1}{\lambda} i\partial\psi\wedge\op\psi,
\end{equation}
{\it the constant in (2.1) could be $\frac{(n-p-q)^2}{8\lambda}$.}

\medskip

\emph{Proof.} By Lemma 7.6, we may assume $u\in C_0^{\infty}(X,\wedge^{n+k}(T_{\mathbb R}^*X\otimes\mathbb C))$, $k\geq 1$. Because
$$
L^k: \ \wedge^{n-k}(T_{\mathbb R}^*X\otimes\mathbb C)\rightarrow \wedge^{n+k}(T_{\mathbb R}^*X\otimes\mathbb C)
$$
is isomorphism, there exists $\varphi \in C_0^{\infty}(X,\wedge^{n-k}(T_{\mathbb R}^*X\otimes\mathbb C))$ such that $u=L^k\varphi$. Thus
$$
u=d\theta+u',
$$
where
$$
\theta=\eta\wedge L^{k-1}\varphi, \ u'=\eta\wedge L^{k-1}d\varphi.
$$
Since
$$
\langle\langle u,u'\rangle\rangle \leq ||u||\sqrt{\lambda} \langle\langle \triangle(L^{k-1}\varphi),L^{k-1}\varphi \rangle\rangle^{1/2},
$$
and
$$
\langle\langle u,d\theta \rangle\rangle \leq ||d^*u||\sqrt{\lambda} ||L^{k-1}\varphi||.
$$
By Lemma 7.5,
$$
\langle\langle \triangle(L^{k-1}\varphi),L^{k-1}\varphi \rangle\rangle \leq \frac 1k \langle\langle \triangle u,u \rangle\rangle, \ ||L^{k-1}\varphi||^2 \leq \frac 1k ||u||^2,
$$
thus
$$
||u||^2=\langle\langle u, d\theta \rangle\rangle-\langle\langle u,u'\rangle\rangle \leq 2\sqrt{\frac{\lambda}{k}}||u||\langle\langle \triangle u,u \rangle\rangle^{1/2}.
$$
(2.1) follows from Lemma 7.3.

\medskip

By Lemma 7.1 and Lemma 7.2,
$$
L^2(X, \wedge^{p,n-p}T^*X)=\mathcal{H}^{p,n-p}_{L^2} \oplus Im(\widetilde{\op}) \oplus  Im((\widetilde{\op})^*),
$$
which shows that every $u\in Dom(\widetilde{\op})\cap Dom(\widetilde{\op^*})$ with degree $p+q = n$, which are orthogonal to the $L^2$ harmonic space, could be written as
$$
u=\widetilde{\op} a+ (\widetilde{\op})^*b, \  a\bot Ker \widetilde{\op}, \ b\bot Ker (\widetilde{\op})^*,
$$
thus
$$
||\widetilde{\op} a||^2\geq\frac{1}{8\lambda}||a||^2, \ ||(\widetilde{\op})^*b||^2\geq\frac{1}{8\lambda}||b||^2,
$$
which proves (2.2).

\medskip

If $\omega$ has a global $C^2$ potential $\psi$, by the Bochner-Kodaira-Nakano identity, for any $u\in C_0^{\infty}(X,\wedge^{p,q}T^*X)$,
\begin{equation}
\int_X(|\op u|^2_{\omega}+|\op^*_{c\psi}u|^2_{\omega})e^{-c\psi} \ dV_{\omega} \geq c(p+q-n)\int_X|u|^2_{\omega}e^{-c\psi} \ dV_{\omega},
\end{equation}
where $c=\frac{1}{2\lambda}(p+q-n)$. Let $u=ve^{\frac c2 \psi}$, then
$$
|\op u|^2_{\omega}e^{-c\psi}=|\op v+\frac c2 \op\psi\wedge v|^2_{\omega}\leq 2|\op v|^2_{\omega}+\frac{\lambda c^2}{2} |v|^2_{\omega},
$$
and
$$
|\op^*_{c\psi}u|^2_{\omega}e^{-c\psi}=|\op^* v+\frac c2 *\partial \psi* v|^2_{\omega}\leq 2|\op^*v|^2_{\omega}+\frac{\lambda c^2}{2} |v|^2_{\omega},
$$
which proves the Theorem finally. $\Box$

\medskip

If $(X, \omega, \lambda)$ satisfies $(2.3)$, we call it \emph{strongly K\"ahler hyperbolic}. In application, we also need the following Theorem. \emph{In the following, we shall use the same symbol $P$ to represent the extension of a differential operator $P$ in the sense of distribution.}

\medskip

\textbf{Theorem 2.2.} {\it Let $(X, \omega, \lambda)$ be a weakly pseudoconvex strongly K\"ahler hyperbolic manifold of dimension $n$, $\varphi$ is plurisubharmonic on $X$. Then for every $v\in L^2_{loc}(X,\wedge^{n,q}T^*X)$, $q \geq 1$, such that $\op v=0$ in the sense of distribution, and}
\begin{equation}
\int_X |v|^2_{\omega}e^{-\varphi} \ dV_{\omega} <+\infty,
\end{equation}
{\it there exist $u\in L^2_{loc}(X,\wedge^{n,q-1}T^*X)$ such that $\op u=v$ in the sense of distribution, and}
\begin{equation}
\int_X |u|^2_{\omega}e^{-\varphi} \ dV_{\omega} \leq \frac{4\lambda}{q^2} \int_X |v|^2_{\omega}e^{-\varphi} \ dV_{\omega}.
\end{equation}

\medskip

\emph{Proof.} The general proof is based on the regularization techniques for plurisubharmonic function (see Demailly \cite{Demailly82}). We will just explain the proof in the simple case when $X$ is bounded pseudoconvex domain in $\mathbb C^n$. By choosing a smooth plurisubharmonic exhaustion function of $X$ and using the standard smoothing technique, we may assume that $\varphi$ and $\psi$ are smooth on some neighborhood of the closure of $X$ in $\mathbb C^n$. Let $u$ be the $L^2$ minimal solution of $\op (\cdot)=v$ in $L^2(X, \wedge^{n,q-1}T^*X, \varphi+\frac{c}{2}\psi)$, where $c=\frac{q}{\lambda}$, then $ue^{\frac c2\psi}$ is the $L^2$ minimal solution of $\op (\cdot)=\op(ue^{\frac c2\psi})$ in $L^2(X, \wedge^{n,q-1}T^*X, \varphi+c\psi)$. Since for every $f\in Dom \ \widetilde{\op}\cap Dom \ (\widetilde{\op})^*_{\varphi+c\psi}$ with degree $(n,q)$,
\begin{equation}
\int_X(|\widetilde{\op} f|^2_{\omega}+|(\widetilde{\op})^*_{\varphi+c\psi}f|^2_{\omega})e^{-(\varphi+c\psi)} \ dV_{\omega} \geq cq \int_X|f|^2_{\omega}e^{-(\varphi+c\psi)} \ dV_{\omega},
\end{equation}
thus
\begin{equation}
\int_X(|ue^{\frac c2\psi}|^2_{\omega}e^{-(\varphi+c\psi)} \ dV_{\omega} \leq \frac{1}{cq} \int_X|v+\frac 2c\op\psi\wedge u|^2_{\omega}e^{-\varphi} \ dV_{\omega} \leq \frac{2}{cq}||v||^2+\frac{c\lambda}{2q}||u||^2,
\end{equation}
we get $(2.6)$. $\Box$

\medskip

By the Theorem 2.1, the lower bound of the Laplace operator depends on $\lambda$. But the Proposition below shows that $\lambda$ has a lower bound.

\medskip

\textbf{Proposition 2.3.} {\it Let $(X,\omega,\lambda)$ be a K\"ahler hyperbolic manifold of dimension $n$, then}
\begin{equation}
\lambda\geq \frac{n}{2n-1}.
\end{equation}

\medskip

\emph{Proof.} By the proof of Lemma 8.4. For any $x\in\widetilde{X}$,
\begin{equation}
    rf'(r)\leq(1+r\sqrt{2n-1}\coth(\frac{r}{\sqrt{2n-1}}))f(r),
\end{equation}
where $f(r)$ is the volume of the geodesic ball $B_x(r)$ in $\widetilde{X}$. We claim that
\begin{equation}
    f'(r)>\sqrt{\frac n\lambda}f(r),
\end{equation}
(which will be proved in the following Lemma). Thus
$$
\frac{r\sqrt{n}}{1+r\sqrt{2n-1}\coth(\frac{r}{\sqrt{2n-1}})}\leq\sqrt{\lambda},
$$
Since $\widetilde{X}$ is complete and non-compact, let $r$ goes to infinity, we get (2.5). $\Box$

\medskip

The following Lemma is contained in Gromov's paper \cite{Gromov91}.

\medskip

\textbf{Lemma 2.4.} {\it Let $(X,\omega)$ be a K\"ahler manifold with the volume form satisfies}
\begin{equation}
\frac{\omega^n}{n!}=d\sigma,
\end{equation}
{\it where $\sigma$ is a $C^1 \ (2n-1)$-form, such that}
\begin{equation}
|\sigma|\leq C,
\end{equation}
{\it on the geodesic ball $B_x(R)$ for some constant $C>0$. If the injectivity radius at point $x$ is bigger than $R$, we have}
\begin{equation}
Cf'(r) \geq f(r),
\end{equation}
{\it for any $r\in [0,R]$.}

\medskip

\emph{Proof.} By using the orthogonal decomposition
\begin{equation}
    \sigma=\sigma_0 ds+\sigma_1
\end{equation}
at every point $z\in\partial B_x(r)$, where $ds$ is the volume form of $\partial B_x(r)$. Thus
$$
|\sigma| \geq |\sigma_0|.
$$
Now
$$
f(r)=\int_{B_x(r)}d\sigma=\int_{\partial B_x(r)}\sigma = \int_{\partial B_x(r)}\sigma_0 ds \leq C f'(r).
$$
which proves the Lemma. $\Box$

\medskip

\textbf{Remark 2.5.} (2.11) follows from
\begin{equation}
|\frac{\eta\wedge\omega^{n-1}}{n!}| \leq \sqrt{\frac{\lambda}{n}}.
\end{equation}
By Bishop's volume comparison theorem, the volume of every compact K\"ahler hyperbolic manifold with Ricci curvature bounded below by some negative constant has a nature upper bound. What we want to show is that Atiyah's $L^2$ index Theorem will give the nature lower bound of $|X|$, which is optimal in some cases.

\medskip

\textbf{Theorem 2.6.} {\it Let $(X,\omega,\lambda)$ be a K\"ahler hyperbolic manifold of dimension $n$, then}
\begin{equation}
    |X| \geq \frac{1}{\sup S_{\widetilde X}}.
\end{equation}

\medskip

\emph{Proof.} By (9.8), $\widetilde{p_1}$ is an integer. By Gromov's result (see Theorem 2.5 in \cite{Gromov91}), $S_{\widetilde X}$ is not identically zero. The Theorem follows from the $\Gamma$ invariance of $S_{\widetilde X}$. $\Box$

\section{Bergman kernels on K\"ahler hyperbolic manifolds}

Denote by $P: \widetilde{X}\to X=\widetilde{X}/\Gamma$ the universal covering map of a complex manifold $X$. According to (7.5) and (7.6), by choosing local coordinate $(z,\bar w)$ of $\widetilde{X} \times \widetilde{X}^*$, one may write
\begin{equation}
B_{\widetilde{X}}^{n,0}=\widetilde K (z,\bar w)dz\otimes d\bar w,
\end{equation}
where $\widetilde K (z,\bar w)$ is locally defined function and $dz,\ d\bar w$ is short for $dz^1 \wedge \cdots \wedge dz^n,\ d\bar w^1 \wedge \cdots \wedge d\bar w^n$ respectively. Using the projection map $P$, $(z,\bar w)$ could also be taken as the local coordinate of $X \times X^*$, i.e., one may write
\begin{equation}
P^*B_{X}^{n,0}=K (z,\bar w)dz\otimes d\bar w,
\end{equation}
(here we use the same notion $P:\ \widetilde{X} \times \widetilde{X}^* \to X\times X^*$ as the canonical projection induced by $P: \widetilde{X}\to X$).

\medskip

Let's recall the notion of tower in Riemannian geometry. Let $\widetilde{X}$ be a Riemannian manifold and $\Gamma$ a free and properly discontinuous group of isometries of $\widetilde{X}$. A tower of subgroups of $\Gamma$ is a nested sequence of subgroups
\begin{equation}
\Gamma=\Gamma_1\supset\Gamma_2\supset\cdots\supset\Gamma_j\supset\cdots\supset\cap\Gamma_j=\{id\},
\end{equation}
such that $\Gamma_j$ is a normal subgroup of $\Gamma$ of finite index $[\Gamma:\Gamma_j]$ for every $j$. The differential manifolds $X_j=\widetilde{X}/\Gamma_j$ are equipped with the push-downs of the Riemannian metric on $\widetilde{X}$. The family $\{X_j\}$ is called a tower of coverings on the Riemannian manifold $X=\widetilde{X}/\Gamma$ (see \cite{DeGeorgeWallach78}). Let
\begin{equation}
    \widetilde{P_j}: \ \widetilde{X}\to X_j, \ P_j: \  X_j\to X, \ P=P_j\circ  \widetilde{P_j}: \ \widetilde{X}\to X,
\end{equation}
be the natural projections. Similar as (9.1) and (9.2), we shall define $F_j(x)$ and $\tau_j(x)$ for very $X_j$. Since $\Gamma_j$ is normal in $\Gamma$, $\tau_j(x)$ must be $\Gamma$ invariant. According to \cite{Donnelly83} (see also \cite{ChenFu} and \cite{DeGeorgeWallach78}), if $X$ is compact,
\begin{equation}
    \tau_j(x)\to \infty,
\end{equation}
uniformly on $\widetilde{X}$ as $j\to \infty$. \emph{Throughout this paper, $\widetilde{X}$ is assumed to be the universal covering of $X$}. It is well known that every Riemannian manifold $X$ with its fundamental group isomorphic to a finitely generated subgroup of $SL(n, \mathbb C)$ admits a tower of coverings with $\widetilde{X}$ being the universal covering (see \cite{Borel63}). Every arithmetic quotient of a bounded symmetric domain satisfies the above condition.

\medskip

Let $X$ be a Hermitian manifold, similar as (3.1), one may write
\begin{equation}
\widetilde{P_j}^*B_{X_j}^{n,0}=\widetilde{K_j}(z,\bar w)dz\otimes d\bar w,\ B_{B_x(\tau_j)}^{n,0}=B_j(z,\bar w)dz\otimes d\bar w,
\end{equation}
and
\begin{equation}
    B_{F_j(x)}^{n,0}=C_j(z,\bar w)dz\otimes d\bar w.
\end{equation}
By (7.9) and (7.10), we have
\begin{equation}
Trace\langle\langle \widetilde{P}_j^*B_{X_j}^{n,0}-B_{F_j(x)}^{n,0}, \widetilde{P}_j^*B_{X_j}^{n,0}-B_{F_j(x)}^{n,0}\rangle\rangle_{F_j(x)}(z)=S_{F_j(x)}(z)-\widetilde{P}_j^*S_{X_j}(z).
\end{equation}
Similarly,
\begin{equation}
Trace\langle\langle B_{\widetilde{X}}^{n,0}-B_{F_j(x)}^{n,0}, B_{\widetilde{X}}^{n,0}-B_{F_j(x)}^{n,0}\rangle\rangle_{F_j(x)}(z) \leq S_{F_j(x)}(z)-S_{\widetilde{X}}(z).
\end{equation}
Therefore
\begin{equation}
Trace\langle\langle B_{\widetilde{X}}^{n,0}-\widetilde{P}_j^*B_{X_j}^{n,0}, B_{\widetilde{X}}^{n,0}-\widetilde{P}_j^*B_{X_j}^{n,0}\rangle\rangle_{F_j(x)} \leq 4(S_{F_j(x)}-S_{\widetilde{X}})
+2(S_{\widetilde{X}}-\widetilde{P}_j^*S_{X_j}).
\end{equation}
Thus we could prove the following Theorem.

\medskip

\textbf{Theorem 1.4.} {\it Let $(X, \omega, \lambda)$ be an n-dimensional non-compact CH K\"ahler hyperbolic manifold with Ricci curvature bounded below by $-1$. Fix $x\in X$, $R>2\sqrt{2n}$, then for every $y$ in the geodesic ball $B_x(R)$ of radius $R$ around $x$, such that $\delta_y:=d(y,\partial B_x(R))\geq 2\sqrt{2n}$, we have}
\begin{equation}
S_{B_x(R)}(y) \leq \sqrt{2}(\frac{2}{\pi})^{n},
\end{equation}
{and}
\begin{equation}
S_{B_x(R)}(y)-S_X(y) \leq 8(\frac{2}{\pi})^{n}\frac{\sqrt{2\lambda}}{\delta_y}.
\end{equation}

\medskip

By a similar argument, one could also get an effective Ramadanov's Theorem for every K\"ahler hyperbolic manifold.

\medskip

\textbf{Theorem 3.1.} {\it Let $\{X_j\}$ be a tower of coverings of a CH K\"ahler hyperbolic manifolds $(X,\omega,\lambda)$ of dimension $n$ and diameter $D$, if $\tau_j > 4(D+2^n\sqrt{n})$, we have}
\begin{equation}
\int_{F_j(x) \times F(x)}|\widetilde{K}(z,\bar w)-K_j(z,\bar w)|^2(i^{n^2}dz\wedge d\bar z)\wedge(i^{n^2}dw\wedge d\bar w) \leq 95|X|(\frac{2}{\pi})^n\frac{\lambda}{\tau_j},
\end{equation}
\begin{equation}
\int_{F_j(x) \times F(x)}|\widetilde{K}(z,\bar w)-C_j(z,\bar w)|^2(i^{n^2}dz\wedge d\bar z)\wedge(i^{n^2}dw\wedge d\bar w) \leq 22|X|(\frac{2}{\pi})^n\frac{\lambda}{\tau_j},
\end{equation}
{\it and}
\begin{equation}
\int_{F_j(x) \times F(x)}|K_j(z,\bar w)-C_j(z,\bar w)|^2(i^{n^2}dz\wedge d\bar z)\wedge(i^{n^2}dw\wedge d\bar w) \leq 26|X|(\frac{2}{\pi})^n\frac{\lambda}{\tau_j}.
\end{equation}

\medskip

\emph{Proof.} By (3.10), the left hand side of (3.11) is no bigger than
\begin{equation}
    I_j:=\int_{F(x)}4(S_{F_j(x)}-S_{\widetilde{X}})
+2(S_{\widetilde{X}}-\widetilde{P}_j^*S_{X_j}).
\end{equation}
By Theorem 1.1 and Theorem 1.4,
\begin{equation}
    I_j \leq 32(\frac{2}{\pi})^n\frac{\sqrt{2\lambda}}{\tau_j-D}|X|
    +32(\frac{4}{\pi})^n\frac{\lambda\sqrt n}{\tau_j^2}|X|.
\end{equation}
By (2.5),
$$
I_j\leq 32(\frac{2}{\pi})^n\frac{\lambda}{\tau_j}|X|(\frac{8}{3}+\frac 1 4)
$$
which proves (3.11). Same method works for (3.12) and (3.13). $\Box$

\medskip

Combining (3.13) with (1.11), we could estimate the Bergman kernels of K\"ahler hyperbolic manifolds.

\medskip

\textbf{Theorem 3.2.} {\it Let $\{X_j\}$ be a tower of coverings of a CH K\"ahler hyperbolic manifolds $(X,\omega,\lambda)$ of dimension $n$ and diameter $D$, if $\tau_j > 4(D+2^n\sqrt{n})$ and $\tau > \sqrt{2n}$, we have}
\begin{equation}
-8(\frac {2}{\pi})^n \frac{\sqrt{2\lambda}}{\tau_j} \leq S_{\widetilde X}-\widetilde{P_j^*}S_{X_j} \leq 26 |X|e^{\frac{\tau^2}{2}}(\frac{2}{\pi^2\tau^2})^{n} \frac{(2n)!}{n!} \frac{\lambda}{\tau_j}.
\end{equation}

\medskip

\emph{Proof.} By (1.11), for every $x\in \widetilde X$,
\begin{equation}
    S_{\widetilde X}(x)-\widetilde{P_j^*}S_{X_j}(x) \geq S_{\widetilde X}(x)-S_{F_j(x)}(x) \geq -8(\frac{2}{\pi})^n\frac{\sqrt{2\lambda}}{\tau_j}.
\end{equation}
By (3.8) and (3.13),
\begin{equation}
    \int_{F(x)}(S_{F_j(x)}-\widetilde{P_j^*}S_{X_j}) \ dV
    \leq 26|X|(\frac{2}{\pi})^n\frac{\lambda}{\tau_j},
\end{equation}
Since the Bergman space of $X_j$ is a closed subspace of the Bergman space of $F_j(x)$, $S_{F_j(x)}-\widetilde{P_j^*}S_{X_j}$ must be a summation of the point-wise norm of some holomorphic n-forms. By Lemma 8.3,
\begin{equation}
   S_{\widetilde X}(x)-\widetilde{P_j^*}S_{X_j}(x) \leq S_{F_j(x)}(x)-\widetilde{P_j^*}S_{X_j}(x)
    \leq 26|X|(\frac{2}{\pi})^{2n}\frac{\lambda}{\tau_j}\frac{e^{\frac{\tau^2}{2}} n!}{(\pi\tau^2)^n} \binom{2n}{n},
\end{equation}
which proves the Theorem. $\Box$

\medskip

One may use the same method to get an effective estimation (depends on $n,\ |B|,\ \lambda, \ D$) of the Bergman kernels on non-CH K\"ahler hyperbolic manifolds. We leave it to the interested reader.

\medskip

The estimation for $K_j(z,\bar w)$ is announced in Yeung \cite{Yeung00b}, we give the details of the proof.

\medskip

\textbf{Theorem 3.3.} {\it Let $\{X_j\}$ be a tower of coverings of a CH K\"ahler hyperbolic manifolds $(X,\omega,\lambda)$ of dimension $n$ and diameter $D$, if $\tau_j > 4(D+2^n\sqrt{n})$, we have}
\begin{equation}
|\widetilde K-K_j|^2(x,y) \leq 190 |X|(\frac{16n}{e\pi^3})^{n} \frac{e^{\frac{\tau^2}{2}}}{\tau^{2n}}  \frac{\lambda}{\tau_j}.
\end{equation}
{\it for every $x, \ y \in X_j$ such that $d(x,y) \leq \frac{\tau_j}{2}+2D$;}
\begin{equation}
|K_j|^2(x,y) \leq 136 |X|(\frac{16n}{e\pi^3})^{n} \frac{e^{\frac{\tau^2}{2}}}{\tau^{2n}}  \frac{\lambda}{\tau_j},
\ \
|\widetilde K|^2(x,y) \leq 128 |X|(\frac{16n}{e\pi^3})^{n} \frac{e^{\frac{\tau^2}{2}}}{\tau^{2n}}  \frac{\lambda}{\tau_j}.
\end{equation}
{\it for every $x, \ y \in X_j$ such that $d(x,y) > \frac{\tau_j}{2}+2D$, where}
\begin{equation}
|K_j|^2(x,y)=|K_j(z,\bar w)|^2 |dz|^2 |dw|^2 (x,y).
\end{equation}

\medskip

\emph{Proof.} If $d(x,y)\leq \frac{\tau_j}{2}+2D$, we have $B_y(\frac{\tau_j}{2}-2D) \subset F_j(x)$, by (3.??)
\begin{equation}
\int_{B_y(\frac{\tau_j}{2}-2D) \times F(x)}|\widetilde{K}-K_j|^2 \ dV \leq 95|X|(\frac{2}{\pi})^n\frac{\lambda}{\tau_j}.
\end{equation}
By (8.3),
$$
\frac{(2\pi n)^n}{\binom{2n}{n}e^n n!} \frac{(\pi \tau^2)^n}{\binom{2n}{n}e^{\frac{\tau^2}{2}} n!}|\widetilde K-K_j|^2(x,y) \leq 95|X|(\frac{2}{\pi})^n\frac{\lambda}{\tau_j}.
$$
Since
$$
\binom{2n}{n} n!=\frac{(2n)!}{n!} \leq 2^{2n+\frac 12}(\frac ne)^n,
$$
we get (3.20). If $d(x,y)>\frac{\tau_j}{2}+2D$. Consider the subdomain
$$
U_{jxy}:=B_x(\frac{\tau_j}{4}+D)\bigcup B_y(\frac{\tau_j}{4}+D)
$$
in $X_j$. Denote by
$$
B_{jxy}^{n,0}:=K_{jxy}(z,\bar w)dz\otimes d\bar w
$$
its Bergman kernel, thus
$$
Trace\langle\langle B_{jxy}^{n,0}-\widetilde{P}_j^*B_{X_j}^{n,0}, B_{jxy}^{n,0}-\widetilde{P}_j^*B_{X_j}^{n,0}\rangle\rangle_{U_{jxy}} \leq S_{U_{jxy}}-\widetilde{P}_j^*S_{X_j}.
$$
Since
$$
B_x(\frac{\tau_j}{4}+D)\bigcap B_y(\frac{\tau_j}{4}+D)=\emptyset,
$$
by definition
$$
K_{jxy}(z,\bar w)=0
$$
for every $z\in B_x(\frac{\tau_j}{4}+D), \ w\in B_y(\frac{\tau_j}{4}+D)$. Thus
$$
((\widetilde{P}_j^*B_{X_j}^{n,0}, \widetilde{P}_j^*B_{X_j}^{n,0}))_{B_y(\frac{\tau_j}{4})\times B_x(\tau)} \leq ((B_{jxy}^{n,0}-\widetilde{P}_j^*B_{X_j}^{n,0}, B_{jxy}^{n,0}-\widetilde{P}_j^*B_{X_j}^{n,0}))_{U_{jxy}\times F(x)}.
$$
By Theorem 1.1 and (1.12), the right hand side is less than
$$
\int_{F(x)}S_{U_{jxy}}-S_{\widetilde X} \ dV + 16|X|(\frac{4}{\pi})^n\frac{\sqrt n \lambda}{\tau_j^2}
\leq 68|X|(\frac{2}{\pi})^n\frac{\lambda}{\tau_j},
$$
which proves the first inequality in (3.21). The second follows by similar method. $\Box$

\section{Very ampleness of the canonical line bundle}

In this section, we will use Bergman metric and Calabi's diastasis function to get an effective very ampleness  criterion the canonical line bundle on K\"ahler hyperbolic manifolds.

\medskip

Let $X$ be a complex manifold, if $X$ has the Bergman metric, i.e.
$$
\mathcal {B}:=i\partial\op \log K(z,\bar z)
$$
is the fundamental form of a K\"ahler metric on $X$, we call $X$ \emph{Bergman hyperbolic}.

\medskip

The notion of \emph{diastasis function} is due to Calabi (see \cite{Calabi53} for definition). The diastasis function for $\mathcal B$ is
\begin{equation}
    \mathcal {D}(z,w):=2\log(\frac{K(z,\bar z)K(w,\bar w)}{|K(z,\bar w)|^2}).
\end{equation}
The notion of \emph{projective-like} is introduced by Loi \cite{Loi06}. We say that $X$ is projective like, if $\mathcal {D}(z,w)|_{\{z \neq w\}}$ has no zero point.

\medskip

The following fundamental Theorem is based on Calabi \cite{Calabi53}.

\medskip

\textbf{Theorem 4.1.} {\it The canonical line bundle $E$ over a compact complex manifold $X$ is very ample if and only if $X$ is Bergman hyperbolic and projective-like.}

\medskip

\emph{Proof.} By definition, $E$ is very ample if and only if all evaluation maps
\begin{equation}
H^0(X,E) \to (J^1 E)_x,  \ H^0(X,E) \to E_x \oplus E_y,\  x, y \in X, \  x \neq y,
\end{equation}
are surjective. Since $Rank(E)=1$, $E$ is very ample if and only if the canonical map $\psi$ from $X$ to $\mathbb P(H^0(X,E))$ is an embedding. Let
$$
\{u_j=u_j(z)dz\}_{j=0,\cdots, N}
$$
be a complete orthonormal base of the Bergman space. The canonical map $\psi$ is given by
\begin{equation}
z \rightarrow [u_0(z),u_1(z),\cdots,u_N(z)],
\end{equation}
where $[u_0,u_1,\cdots,u_N]$ is the homogeneous coordinate of $\mathbb P(H^0(X,E))=\mathbb P^N$. We shall prove the extremal property of the Bergman metric, i.e., for every $Y\in \mathbb C^n \backslash \{0\}$,
\begin{equation}
    \sum_{j,k=1}^n \frac{\partial^2\log K(z,\bar z)}{\partial z_j\partial\bar z_k}Y_j\overline{Y_k}=\frac{1}{K(z,\bar z)}\sup_{||f||=1, \ f(z)=0}\{|Yf(z)|^2 ~|~ f=f(z)dz\in \mathcal H^{n,0}(X)\},
\end{equation}
where
$$
Yf(z)=\sum_{j=1}^n\frac{\partial f}{\partial z_j}(z)Y_j.
$$
Fix $z_0\in X$, one may set $H=\mathcal H^{n,0}(X)$, choose $u_0\in H$ such that $||u_0||=1$ and
$$
|u_0|(z_0)=\sup\{|u|(z_0) ~|~ u\in H, \ ||u||=1\}.
$$
Set $H_1=\{u_0\}^{\bot}\cap H$, choose $u_1^{Y}$ such that $||u_1^Y||=1$ and
$$
|Yu_1^Y(z_0)|=\sup\{~|Yu(z_0)|~ \ |~ u\in H_1, \ ||u||=1\}.
$$
Set $H_2=\{u_0,u_1^Y\}^{\bot}\cap H$. Every $u\in H_1$ satisfies $u(z_0)=0$, every $u\in H_2$ satisfies $Yu(z_0)=0$ and $u(z_0)=0$. Thus by choosing a complete orthornormal base $\{u_2, u_3, \cdots\}$ of $H_2$, one has
\begin{equation}
   K(z_0,\bar z_0) = |u_0(z_0)|^2, \ \sum_{j,k=1}^n \frac{\partial^2\log K(z,\bar z)}{\partial z_j\partial\bar z_k}Y_j\overline{Y_k}=\frac{|Yu_1^Y(z_0)|^2}{|u_0(z_0)|^2},
\end{equation}
which proves the extremal property of the Bergman metric.

\medskip

If $\psi$ is not immersion at $z_0$, the rank of $\psi$ is less than $n$ at $z_0$, thus there exists $Y\in \mathbb C^n\backslash \{0\}$ such that
$$
(Y)_{1\times n}[(\frac{\partial(\frac{u_j}{u_0})}{\partial z_k})(z_0)]_{n\times N}=0,
$$
i.e., $Y\frac{u_j}{u_0}(z_0)=0$ for every $j=1,\cdots, N$. Thus
$$
Yu_1^Y(z_0)=(Y\frac{u_1^Y}{u_0}(z_0))u_0(z_0)=0,
$$
which shows that $X$ is not Bergman hyperbolic.

\medskip

On the other hand, if $X$ is not Bergman hyperbolic, i.e., there exists $z_0\in X$ and $Y\in \mathbb C^n\backslash \{0\}$ such that $Yu_1^Y(z_0)=0$. By definition of $u_1^Y$,
$$
Y\frac{u_j}{u_0}(z_0)=0,
$$
for every $j=1,\cdots, N$. Thus $\psi$ is not immersion at $z_0$. We get that $\psi$ is immersion if and only if $X$ is Bergman hyperbolic.

\medskip

The injection of $\psi$ is equivalent to the following: \\
For every $z,w\in X$, $z\neq w$, $(u_0(z),\cdots, u_N(z))$ is not parallel to $(u_0(w),\cdots, u_N(w))$, i.e.,
$$
(\sum_{j=0}^N|u_j(z)|^2)(\sum_{j=0}^N|u_j(w)|^2)>|\sum_{j=0}^N u_j(z)\overline{u_j(w)}|^2,
$$
which is equivalent to $\mathcal {D}(z,w)>0$. $\Box$

\medskip

By the above Theorem, we could give a very ampleness criterion of the coverings of K\"ahler hyperbolic manifolds, which is a slightly generalization of Yeung's result (see \cite{Yeung00a}).

\medskip

\textbf{Theorem 4.2.} {\it Let $\{X_j\}$ be a tower of coverings of K\"ahler hyperbolic manifold $(X, \omega, \lambda)$. If $\widetilde X$ is weakly pseudoconvex and possess a non-positive strictly plurisubharmonic function $\varphi$. Then the canonical line bundle of $X_j$ is very ample for sufficient large $j$.}

\medskip

\emph{Proof.} \emph{Step 1: $\widetilde X$ is Bergman hyperbolic and for every $\varepsilon>0$, there exists $\delta_{\varepsilon}>0$ such that}
\begin{equation}
   \inf_{p,q\in \widetilde X, \ d(p,q)\geq \varepsilon} \widetilde{\mathcal D}(p,q) \geq \delta_{\varepsilon},
\end{equation}
\emph{where $\widetilde{\mathcal D}$ is the diastatic function of $\widetilde X$.}

\medskip

By Richberg's Theorem (see \cite{Richberg68}), we may assume that $\varphi$ is negative and smooth. Let
$$
\omega_0=i\partial\op(-\log-\varphi).
$$
For every $p\in \widetilde X$, choose a local coordinate $\{z\}$ centered at $p$ such that
$$
\{|z|<1\}\subset\subset\widetilde X.
$$
Take
$$
v_0=\op(\chi(|z|^2)dz),\ v_j=\op(z_j\chi(|z|^2)dz), \ j=1,\cdots,n,
$$
where $\chi\in C^{\infty}(\mathbb R,[0,1])$ is a cut-off function such that $\chi=1$ on $(-\infty,\frac 12)$; $\chi=0$ on $[1,+\infty)$ and $|\chi'|\leq 3$. For every $k\in\mathbb N$, choose sufficient large $C_k$ such that
$$
\varphi_k:=C_k\varphi+2k\chi(|z|^2)\log|z|
$$
is plurisubharmonic on $\widetilde X$. Now
$$
\int_{\widetilde X}|v_0|^2_{\omega_0}e^{-\varphi_n} \ dV_{\omega_0} <\infty, \  \int_{\widetilde X}|v_j|^2_{\omega_0}e^{-\varphi_{n+1}} \ dV_{\omega_0} <\infty, \ j=1,\cdots,n.
$$
Consider the strongly K\"ahler hyperbolic manifold $(\widetilde X,\omega_0)$, by Theorem 2.2, for $j=0,1,\cdots,n$, there exists $u_j\in L^2_{loc}(\widetilde X, \wedge^{n,0}T^*\widetilde X)$ such that $\op u_j=v_j$ in the sense of distribution with $L^2$ estimates. Thus
$$
u_j(0)=0, \ j=0,\cdots,n; \ \frac{\partial u_k}{\partial z_l}(0)=0, \ k,l=1,\cdots,n.
$$
Since $\varphi_k$ is negative and the $L^2$ norm of a $(n,0)$-form not depends on the metric. We get
$$
\chi(|z|^2)dz-u_0, \ z_j\chi(|z|^2)dz-u_j \in \mathcal H^{n,0}(\widetilde X), \ j=1,\cdots,n.
$$
By the extremal property of the Bergman kernel function and the Bergman metric, we get
$$
S_{\widetilde X}(p)>0, \ i\partial\op\log\widetilde K(p)>0,
$$
which proves that $\widetilde X$ is Bergman hyperbolic. By similar method, one could prove that $\widetilde{\mathcal D}(p,q)>0$ for $p \neq q$. Since $S_{\widetilde X}$ is $\Gamma$ invariant, there exists $C>0$ such that
$$
\frac 1C<S_{\widetilde X}<C.
$$
Since
$$
\widetilde{\mathcal D}(p,q)=\widetilde{\mathcal D}(\gamma p,\gamma q)
$$
for every $\gamma\in Aut(\widetilde X)$. In order to prove (4.6), it suffices to show that
\begin{equation}
   \inf_{p\in F(x), \ d(p,q)\geq \varepsilon} \widetilde{\mathcal D}(p,q) \geq \delta_{\varepsilon},
\end{equation}
for fixed $x\in \widetilde X$. By using the same method as in the proof of Theorem 3.3, there exist a sufficient large $G$ such that
$$
|\widetilde K|^2(p,q)\leq\frac{1}{2C^2},
$$
as long as $d(p,q)\geq G$. By definition, if $d(p,q)\geq G$
$$
\widetilde{\mathcal D}(p,q)\geq2\log 2.
$$
Since
$$
U_{G,\varepsilon}:=\{(p,q)\in\overline{F(x)}\times\widetilde X ~|~ \varepsilon\leq d(p,q) \leq G\}
$$
is compact,
$$
\delta_{\varepsilon}:=\min\{2\log2, \inf_{(p,q)\in U_{G,\varepsilon}} \widetilde{\mathcal D}(p,q) \}>0
$$
satisfies (4.7).

\medskip

\emph{Step 2: There exists an constant $A(n,D,|X|,\lambda,\tau)$ such that the canonical line bundle of $X_j$ is very ample for $\tau_j>A$.}

\medskip

By (3.16), there exists $A_1(n,D,|X|,\lambda,\tau)$ such that $X_j$ is Bergman hyperbolic and
$$
\frac 1C<S_{X_j}<C
$$
for $\tau_j>A_1$. Denote by $d_j$ ($\mathcal D_j$) the Bergman distance (Calabi's diastatic) function on $X_j$ respectively. By using the normal coordinate (see Calabi \cite{Calabi53}),
\begin{equation}
    \mathcal D_j(p,q)=d_j(p,q)^2+O(d_j(p,q)^4).
\end{equation}
where $O(d_j(p,q)^4)$ is the curvature terms. Since the Bergman curvature of $X_j$ is bounded (not depends on $j$). There exists $\varepsilon>0$ (not depends on $j$) such that
\begin{equation}
    \mathcal D_j(p,q)\geq\frac 12 d_j(p,q)^2>0
\end{equation}
for $0<d(p,q)\leq\varepsilon$ and $\tau_j>A_1$. Fix such $\varepsilon$, by Theorem 3.3, there exists $A_2$ such that
$$
\frac{S_{X_j}}{S_{\widetilde X}}>(1+\frac{\delta_{\varepsilon}}{4})e^{-\frac{\delta_{\varepsilon}}{2}}, \ |\widetilde K-K_j|^2(p,q)\leq\frac{1}{16C^2}\frac{\delta_{\varepsilon}^2}{e^{\delta_{\varepsilon}}},
$$
for $\tau_j>A_2$. If $d(p,q)\geq\varepsilon$, we have
$$
\frac{\sqrt{S_{X_j}(p)S_{X_j}(q)}}{|K_j|(p,q)}>1,
$$
for $\tau_j>A:=\max\{A_1,A_2\}$, which proves Step 2. $\Box$

\medskip

An elementary method for the effective very ampleness criterion of the compact ball quotients will give in the next section.

\section{Proof of Theorem 1.1}

\emph{Proof of Theorem 1.1.} By Lemma 9.2, we only need to estimate $M^{p,0}_{B_x(\tau)}(x)$, $x\in \widetilde{X}$ and $0\leq p\leq n-1$. Let $f$ be a holomorphic p-forms on $B_x(\tau)$ with
\begin{equation}
\int_{\rho(z)<\tau}|f|^2 \ dV=1,
\end{equation}
where $\rho(z):=d_{\widetilde{X}}(z,x)$. Choose a family of cut-off functions $\chi_{N,\varepsilon}\in C^1(\mathbb R, [0,1])$ such that $\chi_{N,\varepsilon}=1$ on $(-\infty,1-\frac 1N)$ and $\chi_{N,\varepsilon}=0$ on $(1,+\infty)$ with $-N-\varepsilon\leq \chi_{N,\varepsilon}'\leq 0$. By Theorem 2.1, for $0<t<1$,
\begin{eqnarray}
  \nonumber \int_{\rho(z)<\tau(1-t)(1-\frac 1N)}|f|^2 \ dV & \leq & ||\chi_{N,\varepsilon}(\frac{\rho(z)}{\tau(1-t)})f||^2 \\
  \leq \frac{8\lambda}{n-p}||\op(\chi_{N,\varepsilon}(\frac{\rho(z)}{\tau(1-t)})f)||^2 & \leq & \frac{4\lambda(N+\varepsilon)^2}{(n-p)\tau^2(1-t)^2}.
\end{eqnarray}
let $N=2$, $\varepsilon, t$ goes to zero, we have
\begin{equation}
\int_{\rho(z)<\frac{\tau}{2}}|f|^2 \ dV \leq \frac{16\lambda}{(n-p)\tau^2}.
\end{equation}
By the Bochner formula for holomorphic tensor field of covariant degree $p$ (see Kobayashi-Horst \cite{Kobayashi83}),
\begin{equation}
-\op^*\op(|f|^2) \geq -p|f|^2.
\end{equation}
If $n=1$, by Lemma 8.1, Lemma 9.2 and (6.3),
\begin{equation}
\frac{1}{|X|}| \ p_1-\widetilde{p_1}|\leq \frac{64\lambda}{\pi\tau^4}.
\end{equation}
If $n \geq 2$, assume $\tau\geq 2\sqrt{2n}$, by Lemma 8.3 (in case $r=\sqrt{2n}$) and Lemma 9.2,
\begin{equation}
\frac{1}{|X|}| \ p_1-\widetilde{p_1}|\leq
\sum_{p=0}^{n-1}\frac{16\lambda e^n n!}{(n-p)\tau^2 (2\pi n)^n}\binom{n+p}{p} \binom np,
\end{equation}
Since
\begin{equation}
\binom{n+p}{p}=\frac 12 (\binom{n+p}{p} + \binom{n+p}{n}) \leq 2^{n+p-1}, \ \binom np\leq 2^{n-1},
\end{equation}
we get
$$
\frac{1}{|X|}| \ p_1-\widetilde{p_1}|\leq 4(2^n-1) (\frac{2e}{n})^n \frac{n!}{\pi^n}\frac{\lambda}{\tau^2} \leq 16(\frac{4}{\pi})^n\frac{\sqrt{n}\lambda}{\tau^2},
$$
(the second inequality is due to the Stirling's formula) which proves Theorem 1.1. $\Box$

\medskip

By (8.12), Theorem 1.2 follows by using similar method as above.

\section{Proof of Theorem 1.4}

\emph{Proof of Theorem 1.4.} Let $f_y$ be a holomorphic n-form on $B_x(R)$ such that
$$
|f_y|^2(y)=S_{B_x(R)}(y), \ ||f_y||_{B_x(R)}=1.
$$
By Lemma 8.3 (in case $r=\sqrt{2n}$) and the Bochner formula for holomorphic tensor field of covariant degree $n$ (see Kobayashi-Horst \cite{Kobayashi83}), we have
\begin{equation}
    S_{B_x(R)}(y) \leq e^n\frac{n!}{(2\pi n)^n}\binom{2n}{n} \leq (\frac{2}{\pi})^n\sqrt{2},
\end{equation}
(the second inequality is due to the Stirling's formula). Thus
\begin{equation}
    S_X \leq (\frac{2}{\pi})^n\sqrt{2}.
\end{equation}
Choosing a family of cut-off function $\chi_{N,\varepsilon}$ as in the proof of Theorem 1.1. Let $\rho_y(z)=d(y,z)$, for every $0<t<1$, one could solve
\begin{equation}
    \op u_{N,\varepsilon,t,y}=\op(\chi_{N,\varepsilon}(\frac{\rho_y}{(1-t)\delta_y})f_y)
\end{equation}
on $X$ such that
$$
||\chi_{N,\varepsilon}(\frac{\rho_y}{(1-t)\delta_y})f_y||^2=||u_{N,\varepsilon,t,y}||^2+
||\chi_{N,\varepsilon}(\frac{\rho_y}{(1-t)\delta_y})f_y-u_{N,\varepsilon,t,y}||^2,
$$
(i.e., $u_{N,\varepsilon,t,y}$ is the $L^2$ minimal solution. By Theorem 2.1,
\begin{equation}
    ||u_{N,\varepsilon,t,y}||^2 \leq 8\lambda ||\op(\chi_{N,\varepsilon}(\frac{\rho_y}{(1-t)\delta_y})f_y)||^2
    \leq \frac{4\lambda(N+\varepsilon)^2}{(1-t)^2\delta_y^2}.
\end{equation}
Thus
\begin{equation}
    S_X(y) \geq \frac{|\chi_{N,\varepsilon}(\frac{\rho_y}{(1-t)\delta_y})f_y-u_{N,\varepsilon,t,y}|^2(y)}
    {||\chi_{N,\varepsilon}(\frac{\rho_y}{(1-t)\delta_y})f_y-u_{N,\varepsilon,t,y}||^2} \geq S_{B_x(R)}(y)-2\sqrt{S_{B_x(R)}(y)}|u_{N,\varepsilon,t,y}|(y).
\end{equation}
Let $N=2$, $t,\ \varepsilon$ goes to zero, since $u_{N,\varepsilon,t,y}$ is holomorphic on $\rho_y \leq \delta_y (1-t)(1-\frac 1N)$. The Theorem follows by using similar method as in the proof of Theorem 1.1. $\Box$

\section{Appendix I}

We will fix some basic notions of the Hodge theory on non-compact complex manifolds (see Demailly's open-book \cite{Demailly97} and \cite{B02}).

\medskip

Let $(X, \omega)$ be a complex manifold of dimension $n$ and $F_1, \ F_2, \ F_3$ be Hermitian $C^{\infty}$ vector bundles over $X$. Let
$$
P: \\ C^{\infty}(X,F_1)\rightarrow C^{\infty}(X,F_2)
$$
be a differential operator with smooth coefficients, denote by
$$
P^{*}: \\ C^{\infty}(X,F_2)\rightarrow C^{\infty}(X,F_1)
$$
the formal adjoint of $P$, i.e., the unique operator such that for all $u\in C^{\infty}(X,F_1)$ and $v\in C^{\infty}(X,F_2)$,
$$
\langle\langle Pu,v \rangle\rangle = \langle\langle u,P^{*}v \rangle\rangle,
$$
whenever $supp u\cap supp v$ is compact. $P$ induces a non-bounded operator
\begin{equation}
\tilde{P}: \\ L^{2}(X,F_1)\rightarrow L^{2}(X,F_2)
\end{equation}
as follows: we say that $u\in Dom(\tilde{P})$, if there exists a $v(:=\tilde{P}u)\in L^2(X,F_2)$, such that
$$
\langle\langle u,P^{*}g \rangle\rangle = \langle\langle v,g \rangle\rangle
$$
for all $g\in C_{0}^{\infty}(X,F_2)$.

\medskip

$P$ could also induce a non-bounded operator
$$
\vec{P}: \\ L^{2}(X,F_1)\rightarrow L^{2}(X,F_2)
$$
by taking limits. We say that $u\in Dom(\vec{P})$ if there exists a sequence $\{u_j\}\subset C_{0}^{\infty}(X,F_1)$
such that
$$
||u_j-u||\rightarrow 0\\\ (j\rightarrow\infty),
$$
and $\{Pu_j\}$ is a Cauchy sequence of $L^{2}(X,F_2)$. Now
$$
\vec{P}u:=\lim_{j\rightarrow\infty}Pu_j\in L^{2}(X,F_2).
$$
It follows that both $\tilde{P}$ and $\vec{P}$ are densely defined with their graphs closed and
$$
\vec{P}=\tilde{P}|_{Dom\vec{P}}.
$$
By definition,  $C_{0}^{\infty}(X,F_1)$ is dense in $Dom\vec{P}$ for the graph norm with respect to $\vec{P}$. While $C_{0}^{\infty}(X,F_1)$ is dense in $Dom\tilde{P}$ for the graph norm with respect to $\tilde{P}$ if and only if
$\vec{P}=\tilde{P}$. Both $\tilde{P}$ and $\vec{P}$ have a unique Von-Neumann adjoint. Take $\tilde{P}$ for example, its Von-Neumann adjoint
$$
(\tilde{P})^{*}: \\ L^2(X,F_2) \rightarrow L^2(X,F_1)
$$
is defined as follows: we say that $u\in Dom((\tilde{P})^{*})$ if there exists a $v\in L^2(X,F_1)$, such that
$$
\langle\langle u,\tilde{P}f \rangle\rangle = \langle\langle v,f \rangle\rangle
$$
for all $f\in Dom (\tilde{P})$. By definition we have
$$
(\vec{P})^{*}=\overrightarrow{P^{*}}.
$$
But in general, $(\tilde{P})^{*}$ does not coincide with $\widetilde{P^{*}}$. Actually $Dom((\tilde{P})^{*})$ consists of all $u\in Dom(\widetilde{P^{*}})$ satisfying some additional boundary conditions. The following lemma is due to H\"ormander \cite{Hormander65}.

\medskip

\textbf{Lemma 7.1.} {\it The following conditions on $\tilde{P}$ are equivalent:}

\medskip

{\it 1) $Im\tilde{P}$ is closed.}

\medskip

{\it 2) $Im(\tilde{P})^{*}$ is closed.}

\medskip

{\it 3) There is a constant $C$ such that
$$
||u|| \leq C||\tilde{P}u||, \\\ u\in Dom(\tilde{P})\cap\overline{Im(\tilde{P})^{*}}.
$$}

{\it 4) There is a constant $C$ such that
$$
||v|| \leq C||(\tilde{P})^{*}v||, \\\ v\in Dom((\tilde{P})^{*})\cap\overline{Im\tilde{P}}.
$$}
{\it The best constants above are the same.}

\medskip

If there is another differential operator
$$
Q: \\ C^{\infty}(X,F_2) \rightarrow C^{\infty}(X,F_3)
$$
with smooth coefficients, satisfying
$$
\tilde{Q}\circ\tilde{P}=0,
$$
we have
$$
L^2(X,F_2)=H\oplus\overline{Im(\tilde{Q})^{*}}\oplus\overline{Im\tilde{P}}
$$
and
$$
Ker\tilde{Q}=H\oplus\overline{Im\tilde{P}},
$$
where
$$
H=Ker\tilde{Q}\cap Ker(\tilde{P})^{*}.
$$
The following lemma is still due to H\"ormander.

\medskip

\textbf{Lemma 7.2.} {\it A necessary and sufficient condition for $Im(\tilde{Q})^{*}$ and $Im\tilde{P}$ both to be closed is that:}
$$
||u||^2 \leq C^2 (||\tilde{Q}u||^2+||(\tilde{P})^{*}u||^2), \\\ u\in Dom(\tilde{Q}) \cap Dom(\tilde{P})^{*} \cap H^{\bot}.
$$
{\it What's more,}
$$
||u||^2 \leq C^2 (||\tilde{Q}u||^2+||(\tilde{P})^{*}u||^2), \\\ u\in Dom(\tilde{Q}) \cap Dom(\tilde{P})^{*}
$$
{\it is equivalent to}
$$
Im(\tilde{Q})^{*}=Ker(\tilde{P})^{*}, \\\ Im\tilde{P}=Ker\tilde{Q}.
$$

\bigskip

Let $E$ be a rank $r$ holomorphic vector bundle over $X$, with Hermitian metric $h$. Then we have a unique Chern connection $D$ on $\wedge^{p,q}T^{*}X\otimes E$ such that:

\medskip

1) If $D$ is split as a sum of $(1,0)$ and $(0,1)$ connection $D=D'+D''$, $D''=\op$.

\medskip

2) $D$ is Hermitian connection with respect to the metric $h$.

\medskip

Let
$$
\aligned
& \Theta:=D^2:  C^{\infty}(X,\wedge^{p,q}T^{*}X\otimes E) \rightarrow C^{\infty}(X,\wedge^{p+1,q+1}T^{*}X\otimes E),\\
& L:=\omega\wedge\cdot:  C^{\infty}(X,\wedge^{p,q}T^{*}X\otimes E) \rightarrow C^{\infty}(X,\wedge^{p+1,q+1}T^{*}X\otimes E),\\
& \Lambda:=L^{*}:  C^{\infty}(X,\wedge^{p+1,q+1}T^{*}X\otimes E) \rightarrow C^{\infty}(X,\wedge^{p,q}T^{*}X\otimes E),\\
& D:  C^{\infty}(X,\wedge^{r}(T_{\mathbb R}^{*}X\otimes\mathbb C)\otimes E) \rightarrow C^{\infty}(X,\wedge^{r+1}(T_{\mathbb R}^{*}X\otimes\mathbb C)\otimes E),\\
& D^{*}:  C^{\infty}(X,\wedge^{r+1}(T_{\mathbb R}^{*}X\otimes\mathbb C)\otimes E) \rightarrow C^{\infty}(X,\wedge^{r}(T_{\mathbb R}^{*}X\otimes\mathbb C)\otimes E),\\
& \op:  C^{\infty}(X,\wedge^{p,q}T^{*}X\otimes E) \rightarrow C^{\infty}(X,\wedge^{p,q+1}T^{*}X\otimes E),\\
& \op^{*}:  C^{\infty}(X,\wedge^{p,q+1}T^{*}X\otimes E) \rightarrow C^{\infty}(X,\wedge^{p,q}T^{*}X\otimes E),\\
& D':  C^{\infty}(X,\wedge^{p,q}T^{*}X\otimes E) \rightarrow C^{\infty}(X,\wedge^{p+1,q}T^{*}X\otimes E),\\
& (D')^{*}:  C^{\infty}(X,\wedge^{p+1,q}T^{*}X\otimes E) \rightarrow C^{\infty}(X,\wedge^{p,q}T^{*}X\otimes E),\\
& \Delta':=D'(D')^{*}+(D')^{*}D':  C^{\infty}(X,\wedge^{p,q}T^{*}X\otimes E) \rightarrow C^{\infty}(X,\wedge^{p,q}T^{*}X\otimes E),\\
& \Delta'':=\op\op^{*}+\op^{*}\op:  C^{\infty}(X,\wedge^{p,q}T^{*}X\otimes E) \rightarrow C^{\infty}(X,\wedge^{p,q}T^{*}X\otimes E),\\
& \Delta:=DD^{*}+D^{*}D: C^{\infty}(X,\wedge^{r}(T_{\mathbb R}^{*}X\otimes\mathbb C)\otimes E) \rightarrow C^{\infty}(X,\wedge^{r}(T_{\mathbb R}^{*}X\otimes\mathbb C)\otimes E).
\endaligned
$$

\medskip

The Hermitian metric $h$ on $E$ could induce a sesquilinear pairing $\{\bullet,\bullet\}$ as follows. For an arbitrary holomorphic trivialization
$$
\theta: \ E|_{\Omega}\rightarrow \Omega\times\mathbb C^r.
$$
Let $H=(h_{\lambda\overline{\mu}})$ be the Hermitian matrix with smooth coefficients representing the metric along the fibre of $E|_{\Omega}$. For any $s,t \in C^{\infty}_{\bullet,\bullet}(X,E)$ and $\sigma=\theta(s), \tau=\theta(t)$, one can write
$$
\{s,t\}=\sum_{\lambda,\mu}h_{\lambda\overline{\mu}} \sigma^{\lambda} \wedge \overline{\tau}^{\mu}.
$$
The Hodge-Poincar\'e-De Rham operator $*$ is the collection of $\mathbb{C}$-linear isometric maps defined by
$$
* \ : \wedge^{p,q}T^{*}X\otimes E \rightarrow \wedge^{n-q,n-p}T^{*}X\otimes E, \  \{s,*t\}=\langle s,t\rangle dV,
$$
where $\langle\bullet,\bullet\rangle$ is the point-wise Hermitian inner product.
The following lemma is classical:

\medskip

\textbf{Lemma 7.3.} {\it $*\triangle''=\triangle'*, \ *\triangle'=\triangle''*, \ *\triangle=\triangle*$. What's more, if $(X,\omega)$ is K\"ahler and $E$ is trivial with trivial metric, $\triangle=2\triangle'=2\triangle''$ commutes with all operators $*,\partial,\partial^*,\op,\op^*,L,\Lambda$.}

\medskip

Since the order of $L,\Lambda$ is zero, they could be point-wisely defined. Let's recall the notion of primitive: a homogeneous element $u\in\wedge^s(T_{\mathbb R}^*X\otimes\mathbb C)$ is called primitive if $\Lambda u=0$. The space of primitive elements of total degree s will be denoted by
$$
Prim^sT^*X=\oplus_{p+q=s}Prim^{p,q}T^*X.
$$
Now we can state the primitive decomposition formula (see \cite{Demailly97} for the proof):

\medskip

\textbf{Lemma 7.4.} {\it For every $u\in \wedge^s(T_{\mathbb R}^*X\otimes\mathbb C)$, there is a unique decomposition}
$$
u=\sum_{r=(s-n)_+}^{[s/2]}L^r u^{(s-2r)},  \  u^{(s-2r)}\in Prim^{s-2r}T^*X,
$$
{\it where $(s-n)_+=\max\{0,s-n\}$, $[s/2]$ is the integer part of $s/2$.}

\medskip

With the help of primitive decomposition formula, we could prove the following result which is crucial in finding the precise lower bound of the spectrum of Laplace-Beltrami operator in K\"ahler hyperbolic case.

\medskip

\textbf{Lemma 7.5.} {\it If $(X,\omega)$ is K\"ahler, $\varphi\in C_0^{\infty}(X, \wedge^{n-k}(T_{\mathbb R}^*X\otimes\mathbb C))$, $k \geq 1$.}
$$
\langle\langle \triangle(L^k \varphi),L^k \varphi\rangle\rangle\geq k \langle\langle \triangle(L^{k-1} \varphi),L^{k-1} \varphi\rangle\rangle, \  ||L^k \varphi||^2 \geq k ||L^{k-1} \varphi||^2.
$$

\medskip

\emph{Proof.} For any $u\in C_0^{\infty}(X, \wedge^{s}(T_{\mathbb R}^*X\otimes\mathbb C))$, by primitive decomposition formula,
$$
u=\sum_{r=(s-n)_+}^{[s/2]}L^r u^{(s-2r)}.
$$
Since $u^{(s-2r)}$ is primitive, according to Lemma 7.3,
$$
\langle\langle \triangle u, u\rangle\rangle=\sum_{r=(s-n)_+}^{[s/2]}(\prod_{t=1}^r t(n+2r-s-t+1))\langle\langle \triangle u^{(s-2r)}, u^{(s-2r)}\rangle\rangle,
$$
and
$$
||u||^2=\sum_{r=(s-n)_+}^{[s/2]}(\prod_{t=1}^r t(n+2r-s-t+1))||u^{(s-2r)}||^2,
$$
where $\prod_{t=1}^0:=1$. Thus for $0\leq m\leq k$,
$$
\langle\langle \triangle(L^m \varphi),L^m \varphi\rangle\rangle=\sum_{r=0}^{[n-k/2]} b_{k,m,r}
\langle\langle \triangle\varphi^{(n-k-2r)},\varphi^{(n-k-2r)}\rangle\rangle,
$$
and
$$
||L^m \varphi||^2=\sum_{r=0}^{[n-k/2]} b_{k,m,r}||\varphi^{(n-k-2r)}||^2,
$$
where
$$
b_{k,m,r}=\prod_{t=1}^{r+m}t(2r+k+1-t).
$$
Hence
$$
||L^k \varphi||^2 \geq \min_{0\leq r\leq [n-k/2]}\{\frac{b_{k,k,r}}{b_{k,k-1,r}}\} ||L^{k-1} \varphi||^2=k ||L^{k-1} \varphi||^2,
$$
and
$$
\langle\langle \triangle(L^k \varphi),L^k \varphi\rangle\rangle\geq k \langle\langle \triangle(L^{k-1} \varphi),L^{k-1} \varphi\rangle\rangle.
$$
The Lemma is proved. $\Box$

\medskip

By the above discussion, one can define $\tilde{\Theta},\tilde{L},\cdots$ as closed and densely defined operators on $L^2_{\bullet}:=\oplus_{r=1}^{2n}L^2(X,\wedge^{r}(T_{\mathbb R}^{*}X\otimes\mathbb C)\otimes E)
=\oplus_{p,q=1}^{n}L^2(X,\wedge^{p,q}T^{*}X\otimes E)$. Because $\op\op=0$, by definition, $\widetilde{\op}\widetilde{\op}=0$, one has the orthogonal decompositions
\begin{equation}
L^2_{\bullet}=\mathscr{H}^{\bullet}_{\widetilde{\op}} \oplus \overline{Im(\widetilde{\op})} \oplus \overline{Im((\widetilde{\op})^*)},
\end{equation}
where
$$
\mathscr{H}^{\bullet}_{\widetilde{\op}}=Ker(\widetilde{\op}) \cap Ker((\widetilde{\op})^*) \subset C^{\infty}_{\bullet},
$$
(Similar results and definitions for $\mathscr{H}^{\bullet}_{\widetilde{\op^*}}$, $\mathscr{H}^{\bullet}_{\widetilde{D'}}$ and $\mathscr{H}^{\bullet}_{\widetilde{(D')^{*}}}$). We also define
$$
\mathscr{H}^{\bullet}_{\widetilde{D}}=Ker(\widetilde{D}) \cap Ker((\widetilde{D})^*) \subset C^{\infty}_{\bullet}.
$$
If $E$ is flat, one has also
\begin{equation}
L^2_{\bullet}=\mathscr{H}^{\bullet}_{\widetilde{D}} \oplus \overline{Im(\widetilde{D})} \oplus \overline{Im((\widetilde{D})^*)}.
\end{equation}
The space of $L^2$ harmonic forms with respect to $\Delta$ will be defined as $\mathcal {H}^{\bullet}_{\widetilde{\Delta}}=Ker(\widetilde{\Delta})$, (similar definitions for $\Delta'$ and $\Delta''$). If $(X,\omega)$ is non-complete, generally $\mathscr{H}^{\bullet}_{\widetilde{D}}$ is only a subset of $\mathcal {H}^{\bullet}_{\widetilde{\Delta}}$, but if $(X,\omega)$ is complete, one has the following classical result (see H\"ormander's density technique \cite{Hormander65}, \cite{AndreottiVesentini65} and Demailly's open book \cite{Demailly97}).

\medskip

\textbf{Lemma 7.6.} {\it If $(X,\omega)$ is complete, one has}
\begin{equation}
\widetilde{D}=\vec{D}, \ Dom\widetilde{\Delta} \subset Dom(\widetilde{D}) \cap Dom(\widetilde{D^*}),
\end{equation}
{\it (Similar results for $\triangle'$ and $\triangle''$) and}
\begin{equation}
\mathscr{H}^{\bullet}_{\widetilde{D}}=\mathscr{H}^{\bullet}_{\widetilde{D^*}}=\mathcal{H}^{\bullet}_{\widetilde{\Delta}}:
=\mathcal{H}^{\bullet},
\end{equation}
{\it (Similar definitions for $\mathcal{H}^{\bullet}_{'}$ and $\mathcal{H}^{\bullet}_{''}$).}

\medskip

The Schwartz kernel for the projection onto $\mathcal{H}^{p,q}_{''}(X,E)$ is defined as
\begin{equation}
B_{X,E}^{p,q}:=\sum_j u_j\boxtimes\overline{u_j}\in H^0(X\times X^*,(\wedge^{p,q}T^*X\otimes E)\boxtimes(\wedge^{p,q}T^*(X^*)\otimes \overline E)):=H^{p,q}_{\boxtimes},
\end{equation}
where $u_j$ is any complete orthonormal base of the separated Hilbert space $\mathcal{H}^{p,q}_{''}(X,E)$, $X^*$ the conjugated complex manifold of $X$ and $\overline{E}$ the conjugate of $E$. Now let's recall the notion $\boxtimes$. Let $L_X$ and $L_Y$ be the vector bundles over the complex manifolds $X$ and $Y$ respectively. By definition,
\begin{equation}
L_X \boxtimes L_Y:=P_1^* L_X \otimes P_2^* L_Y,
\end{equation}
where $P_1: \ X \times Y\to X$ and $P_2: \ X \times Y\to Y$ are the canonical projection map.

\medskip

For any $u\in \mathcal{H}^{p,q}_{''}(X,E)$ and $P,Q,R\in H^{p,q}_{\boxtimes}$, we shall define
\begin{equation}
\langle\langle u,P \rangle\rangle=\sum_{j,k}\langle\langle u(\cdot), P^{\bar{k}j}(\cdot,\bar w)P_j(\cdot)\rangle\rangle P_k(w),
\end{equation}
\begin{equation}
\langle\langle P,Q \rangle\rangle=\sum_{j,k,l,m}\langle\langle P^{\bar{k}j}(\cdot,\bar w)P_j(\cdot), Q^{\bar{m}l}(\cdot,\bar z)P_l(\cdot)\rangle\rangle P_m(z)\boxtimes\overline{P_k(w)},
\end{equation}
\begin{equation}
Trace(R)=\sum_{j,k}R^{\bar kj}(z,\bar z)\langle R_j(z), R_k(z)\rangle,
\end{equation}
and
\begin{equation}
((P,Q))=\int_X Trace\langle\langle P,Q \rangle\rangle dV,
\end{equation}
where
$$
P=\sum_{j,k}P^{\bar kj}(z,\bar w) P_j(z)\boxtimes \overline{P_k(w)},
$$
(Similar formula for $Q$ and $R$). Then we have the reproducing formula,
\begin{equation}
\langle\langle u,B^{p,q}_{X,E} \rangle\rangle=u.
\end{equation}
Thus, sometimes we call $B^{p,q}_{X,E}$ the reproducing kernel or generalized Bergman kernel of $\mathcal{H}^{p,q}_{''}(X,E)$.

\medskip

We call
\begin{equation}
S_{X,E}^{p,q}:=Trace(B^{p,q}_{X,E})
\end{equation}
and
\begin{equation}
\mathcal {S}_{X,E}^{p,q}:=Trace(B^{p,q}_{X,E})dV
\end{equation}
the Schwartz kernel functions and forms respectively. When $E$ is trivial with trivial metric, we shall omit the lower index $E$, and call
\begin{equation}
S_X:=S_X^{n,0}
\end{equation}
and
\begin{equation}
\mathcal {S}_X:=\mathcal {S}_X^{n,0}
\end{equation}
the Bergman kernel function and form respectively.

\medskip

Sometimes, it is convenient to use the following extremal function:
\begin{equation}
M_{X,E}^{p,q}:=\sup\{|u|^2 \ | \ u\in \mathcal{H}^{p,q}_{''}(X,E), \ ||u||=1\},
\end{equation}
(similar definition for $M_X^{p,q}$). By Berndtsson's Lemma \cite{Berndtsson02},
\begin{equation}
M_{X,E}^{p,q}\leq S_{X,E}^{p,q}\leq r \binom{n}{p}\binom{n}{q}M_{X,E}^{p,q},
\end{equation}
where $\binom{n}{p}$ are binomial coefficients.

\medskip

If $(X,\omega)$ is a K\"{a}hler manifold, one has the Bochner-Kodaira-Nakano identity (see \cite{OhsawaTakegoshi87} for generalization)
\begin{equation}
\Delta''=\Delta'+[i\Theta,\Lambda],
\end{equation}
which implies:

\medskip

\textbf{Lemma 7.7.} {\it If $(X,\omega)$ is a complete K\"{a}hler manifold and $E$ is flat, one has the following orthogonal decomposition}
\begin{equation}
\mathcal {H}^{r}_{L^2}(M,E)=\oplus_{p+q=r}\mathcal {H}^{p,q}_{L^2}(M,E),
\end{equation}
where
\begin{equation}
\mathcal{H}^{\bullet}=\mathcal{H}^{\bullet}_{'}=\mathcal{H}^{\bullet}_{''}:=\mathcal{H}^{\bullet}_{L^2}.
\end{equation}

\section{Appendix II}

A Cartan-Hadamard (CH) manifold is a complete, simply-connected Riemannian manifold of nonpositive curvature. The following result is due to Green and Wu \cite{GreeneWu79}.

\medskip

\textbf{Lemma 8.1.} {\it Let $X$ be a CH manifold with (real) dimension $2n$, for any point $x$ in $X$, if $u$ is a non-negative subharmonic function on the geodesic open ball $B_x(r)$ of radius $r$ around $x$ in $X$, then}
\begin{equation}
\int_{B_x(r)}u \ dV \geq \frac{(\pi r^2)^n}{n!}u(x).
\end{equation}

\medskip

Another version is due to Li-Schoen \cite{L-S84}.  Theorem 2.1 in \cite{L-S84} is slightly different from the following Theorem. But the Theorem below follows easily from Li-Schoen's result by estimating the constant in Theorem 2.1 \cite{L-S84} carefully.

\medskip

\textbf{Lemma 8.2.} {\it Let $X$ be a $m$ real dimensional compact Riemannian manifold with (possibly empty) boundary, for any point $x$ in $X$, if $X$ has no boundary, assume its diameter is no less than $2r$, otherwise assume that the distance from $x$ to the boundary of $X$ is at least $5r$. Suppose the Ricci curvature of $X$ is bounded below by $-(m-1)k^2$, where $k\geq 0$. If $u$ is a non-negative subharmonic function on the geodesic open ball $B_x(r)$ of radius $r$ around $x$ in $X$, then}
\begin{equation}
\int_{B_x(r)}u \ dV \geq 2^{-200(m-1)(1+kr)}|B_x(r)|\sup_{B_x(r/2)}u,
\end{equation}
where $|B_x(r)|$ stands for the volume of $B_x(r)$.

\medskip

We also have similar results on non-subharmonic functions.

\medskip

\textbf{Lemma 8.3.} {\it Let $X$ be a CH K\"ahler manifold with (complex) dimension $n$. For any non-negative smooth function $u$ on the geodesic open ball $B_x(r)$ of radius $r$ around $x$ in $X$, if}
$$
-\op^*\op u+pu\geq 0,
$$
{\it on $B_x(r)$ for some $p\in \mathbb Z^+$, then}
\begin{equation}
u(x)\leq e^{r^2/2}\frac{n!}{(\pi r^2)^n}\binom{n+p}{p}\int_{B_x(r)}u \ dV.
\end{equation}

\medskip

\emph{Proof.} Consider the complete K\"ahler metric
\begin{equation}
    2Re(\sum_{j=1}^{p}dt^j\otimes d\bar{t}^j+\sum_{\alpha,\beta=1}^{n}g_{\alpha\overline{\beta}}dz^\alpha\otimes d\bar{z}^\beta),
\end{equation}
on $\mathbb C^p\times X$, where $\omega=i\sum_{\alpha,\beta=1}^{n}g_{\alpha\overline{\beta}}dz^\alpha\wedge d\bar{z}^\beta$. Denote by $\omega_t$, $dV_t$ the associated fundamental form and the volume form respectively. Now our distance function
\begin{equation}
    d_t((0,x),(t,z))=\sqrt{2|t|^2+\rho(z)^2},
\end{equation}
where $\rho(z)=d(x,z)$ is the distance function on $(X,\omega)$. Denote by $\op_t$ the Cauchy-Riemann operator on $\mathbb C^p\times X$, we have
\begin{equation}
    -\op_t^*\op_t(e^{|t|^2}u) \geq |t|^2e^{|t|^2}u\geq 0,
\end{equation}
by Lemma 8.1, we have
\begin{equation}
   \int_{d_t<r}e^{|t|^2}u \ dV_t \geq \frac{(\pi r^2)^{n+p}}{(n+p)!}u(x).
\end{equation}
The Lemma follows from
$$
\int_{d_t<r}e^{|t|^2}u \ dV_t \leq \int_{\{2|t|^2<r^2\}\times B_x(r)}e^{|t|^2}u \ dV_t
$$
$$
\leq (2^p e^{r^2/2}\int_{2|t|^2<r^2}d\lambda)(\int_{B_x(r)}u \ dV)=e^{r^2/2}\frac{(\pi r^2)^{p}}{p!}\int_{B_x(r)}u \ dV,
$$
where $d\lambda$ is the Lebesgue measure on $\mathbb C^p$. $\Box$

\medskip

One needs the following Lemma to generalize Lemma 8.2.

\medskip

\textbf{Lemma 8.4.} {\it Let $X$ be a $m$ real dimensional compact (without boundary) Riemannian manifold with Ricci curvature bounded below by $-(m-1)k^2$, $k>0$. Denote by $D$ the diameter of $X$. For any $x\in X$, if $0<a \leq b \leq D$ and $\partial B_x(b)\neq\emptyset$, we have}
\begin{equation}
\frac{|B_x(b)|}{|B_x(a)|} \leq (\frac{b}{a})(\frac{\sinh(bk)}{\sinh(ak)})^{m-1} \leq (\frac{b}{a})^m e^{(m-1)(b-a)k}.
\end{equation}

\medskip

\emph{Proof.} Let $\rho(\cdot):=d(\cdot,x)$ be the distance function of $X$. By Laplacian comparison Theorem (see Green and Wu \cite{GreeneWu79}),
\begin{equation}
    \triangle \rho\leq (m-1)k\coth(k\rho)
\end{equation}
in the sense of distribution. Thus
\begin{equation}
\triangle \rho^2 \leq 2+2(m-1)k\rho\coth(k\rho).
\end{equation}
For any $0<r\leq b$, integrating the above formula over $B_x(r)$, we get
$$
\int_{B_x(r)}\triangle \rho^2 \ dV \leq (2+2(m-1)kr\coth(kr))|B_x(r)|.
$$
Since
\begin{equation}
\int_{B_x(r)}\triangle \rho^2 \ dV=\int_{\partial B_x(r)} * d\rho^2 = 2r(\frac{d}{dr}|B_x(r)|),
\end{equation}
we have
$$
\frac{|B_x(b)|}{|B_x(a)|} \leq e^{\int_a^b(1/r+(m-1)k\coth(kr))dr} = (\frac{b}{a})(\frac{\sinh(bk)}{\sinh(ak)})^{m-1},
$$
which proves the Lemma. $\Box$

\medskip

Consider the product $X\times\mathbb R$, by the above Lemma, we have the similar result as Lemma 8.3.

\medskip

\textbf{Lemma 8.5.} {\it Let $(X, \omega)$ be a compact K\"ahler manifold of complex dimension n, with Ricci curvature bounded below by $-1$ and diameter no less than $2$. For any non-negative smooth function $u$ on the unit ball $B_x$ with center $x$ in $X$, if}
$$
-\op^*\op u+pu\geq 0,
$$
{\it on $B_x$ for $0 \leq p \leq n-1$, then}
\begin{equation}
u(x)\leq 2^{201(2n+\sqrt{2n})}\frac{1}{|B_x|}\int_{B_x}u \ dV.
\end{equation}

\medskip

\section{Appendix III}

Let's recall some basic facts on covering spaces. Let $(\widetilde{X}, \tilde{g})$ be a Riemannian manifold. Let $\Gamma$ be a subgroup of the isometrics that acts freely and properly discontinuously on $\widetilde{X}$. Let $X=\widetilde{X}/\Gamma$ be the quotient manifold and $P: \widetilde{X}\rightarrow X$ be the covering map. Equip $X$ with the push-down metric $g$ so that $P^*g=\tilde{g}$. Denote by $d_{\widetilde{X}}$ and $d_X$ the distance functions of $(\widetilde{X},\tilde{g})$ and $(X,g)$ respectively. For $x\in \widetilde{X}$, let
\begin{equation}
    F(x):=\{y\in\widetilde{X} \ | \ d_{\widetilde{X}}(y,x)<d_{\widetilde{X}}(y,\gamma x), \forall \gamma \in \Gamma\setminus 1\}
\end{equation}
be the Dirichlet fundamental domain centered at $x$ and
\begin{equation}
    \tau(x):=\frac12\inf\{d_{\widetilde{X}}(x,\gamma x) \ | \ \gamma \in \Gamma\setminus 1 \}
\end{equation}
be the quasi-injectivity radius of $P(x)$ in $X$ with respect to the covering map $P$. By definition, the geodesic ball $B_x(\tau(x))$ is contained in $F(x)$ and if $\widetilde{X}$ has no conjugate points, $\tau(x)$ is the injectivity radius of $P(x)$ in $X$. In particular, this is the case when $\widetilde{X}$ is CH manifold.

\medskip

If $(\widetilde{X},\widetilde{\omega})$ is the universal covering space of  the compact complex manifold $(X, \omega)$, one has the following famous $L^2$ index formula due to Atiyah \cite{Atiyah76}.

\medskip

\textbf{Lemma 9.1.} {\it For any compact complex manifold $(X, \omega)$, one has}
\begin{equation}
\chi(X,\Omega^p)=\chi_{L^2}(X,\Omega^p).
\end{equation}

\medskip

Here $\Omega^p$ is the sheaf of germs of holomorphic p-forms on $n$ dimensional compact complex manifold $(X, \omega)$ . The holomorphic Euler characteristic of $\Omega^p$  is defined as
\begin{equation}
    \chi(X,\Omega^p):=\sum_{q=0}^n (-1)^q h^{p,q},
\end{equation}
where the Hodge number $h^{p,q}$ is defined as
$$
\int_X \mathcal{S}_X^{p,q}.
$$
The $L^2$ Euler characteristic of $\Omega^p$  is defined as
\begin{equation}
    \chi_{L^2}(X,\Omega^p):=\sum_{q=0}^n (-1)^q \widetilde{h^{p,q}},
\end{equation}
where the $L^2$ Hodge number $\widetilde{h^{p,q}}$ is defined as
$$
\int_X \mathcal{S}_{\widetilde X}^{p,q},
$$
where $(\widetilde X, \widetilde{\omega}:=P^*\omega)$ is the universal covering complex manifold of $(X, \omega)$.

\medskip

According to Atiyah's $L^2$ Index Formula, if $X$ is K\"ahler hyperbolic, we have the following Lemma.

\medskip

\textbf{Lemma 9.2.} {\it Let $X$ be a K\"ahler hyperbolic manifold, one has}
\begin{equation}
    \frac{1}{|X|}| \ p_1-\widetilde{p_1}|\leq \sum_{p=0}^{n-1}\binom np \sup_{z\in\widetilde X}M_{B_z(\tau)}^{p,0}(z).
\end{equation}

\medskip

\emph{Proof.} According to Atiyah's $L^2$ Index Formula,
\begin{equation}
\sum_{p=0}^n (-1)^p h^{0,p}=\sum_{p=0}^n (-1)^p \widetilde{h^{0,p}}.
\end{equation}
By Lemma 7.7 and Lemma 2.1,
\begin{equation}
\widetilde{p_1}-p_1=\sum_{p=0}^{n-1} (-1)^{n+p} h^{p,0}.
\end{equation}
For any $x\in \widetilde X$, by definition
$$
h^{p,0}=\int_{F(x)}P^*S_X^{p,0} \ dV,
$$
by Berndtsson's lemma, for any $z\in \widetilde X$,
$$
P^*S_X^{p,0}\leq\binom np P^*M_X^{p,0}(z).
$$
And by definition, one has
$$
P^*M_X^{p,0}(z)\leq P^*M_{P(F(z))}^{p,0}(z)=M_{F(z)}^{p,0}(z),
$$
and
$$
M_{F(z)}^{p,0}(z)\leq M_{B_z(\tau(z))}^{p,0}(z)\leq M_{B_z(\tau)}^{p,0}(z).
$$
So the Lemma is proved. $\Box$

\section*{Acknowledgment}

The results of this paper were obtained during my Ph.D. studies at Tongji University. I would like to express deep gratitude to my supervisor Bo-Yong Chen whose guidance and support were crucial for the completion of this paper.


\begin{thebibliography}{99}

\bibitem{AndreottiVesentini65} A. Andreotti and E. Vesentini, {\it Carleman estimates for the
    Laplace-Beltrami equation in complex manifolds}, Publ. Math. IHES {\bf 25} (1965), 81--130.

\bibitem{Atiyah76} M. F. Atiyah, {\it Elliptic operators, discrete groups and von Neumannn algebras},
    Ast\'erisque {\bf 32-33} (1976), 43--72.

\bibitem{Berndtsson02} B. Berndtsson, {\it An eigenvalue estimate for the $\bar{\partial}-$Laplacian}, J.
    Diff. Geom. {\bf 60} (2002), 295--313.

\bibitem{B02} J.~Bertin, J.-P. Demailly, L.~Illusie, and C.~Peters, \emph{Introduction to Hodge theory}, SMF/AMS Texts and Monographs, vol.~8, American Mathematical Society, Providence, RI, 2002.

\bibitem{Borel63} A. Borel, {\it Compact Clifford-Klein forms of symmetric spaces}, Topology {\bf 2}(1963),
    111-122.

\bibitem{Calabi53} E. Calabi, {\it Isometric imbeddings of complex manifolds}, Ann. of Math. {\bf 58} (1953), 1--23.

\bibitem{ChenFu} B.-Y. Chen and S. Fu, {\it Stability of the Bergman kernel on a tower of coverings}, preprint.

\bibitem{DM83} J. E. D'atri and I. Dotti. Miatello, {\it A characterization of bounded symmetric domains by curvature}, Trans. Am. Math. Soc. {\bf 276} (1983), 531--540.

\bibitem{DeGeorgeWallach78} D. DeGeorge and N. Wallach, {\it Limit formulas for multiplicities in
    $L^2(\Gamma\backslash G)$}, Ann. of Math. {\bf 107} (1978), 133--150.

\bibitem{Demailly82} J.-P. Demailly, {\it Estimations $L^2$ pour l'op\'{e}rateur $\bar{\partial}$ d'un
    fibr\'{e} vectoriel holomorphe semi-positif au-dessus d'une vari\'{e}t\'{e} k\"{a}hl\'{e}rienne
    compl\`{e}te}, Ann. Sci. \'{E}cole Norm. Sup. {\bf 15} (1982), 457--511.

\bibitem{Demailly97} J.-P. Demailly, {\it Complex analytic and differential geometry}. Book available from the author's homepage.

\bibitem{Donnelly83} H. Donnelly, {\it On the spectrum of towers}, Proc. Amer. Math. Soc. {\bf 87} (1983),
    322--329.

\bibitem{Donnelly96} H. Donnelly, {\it Elliptic operators and covering of Riemmannian manifolds}, Math. Z.
    {\bf 233} (1996), 303--308.

\bibitem{DonnellyFefferman83} H. Donnelly and C. Fefferman, {\it $L^2-$cohomology and index theorem for the
    Bergman metric}, Ann. of Math. {\bf 118} (1983), 593--618.

\bibitem{GreeneWu79} R. E. Greene and H. Wu, Function Theory on Manifolds Which Possess a Pole, Lect. Notes
    in Math. {\bf 699}, 1979.

\bibitem{Gromov91} M. Gromov, {\it K\"ahler hyperbolicity and $L^2$-Hodge theory}, J. Diff. Geom. {\bf 33}
    (1991), 263--292.

\bibitem{H-P96} S. Hersonsky and F. Paulin, {\it On the volumes of complex hyperbolic manifolds}, Duke. Math. J. {\bf 84} (1996), 719--737.

\bibitem{Hormander65} L. H\"ormander, {\it $L^2$ estimates and existence theorems for the $\op$ operators}, Acta. Math. {\bf 113} (1965), 89--152.

\bibitem{Ishi} H. Ishi, {\it On the Bergman metric of Siegel domains}, preprint.

\bibitem{Kai_Ohsawa07} C. Kai and T. Ohsawa, {\it A note on the Bergman metric of the bounded homogeneous domains}, Nagoya Math. J. {\bf 186}
    (2007), 157--163.

\bibitem{Kazhdan70} D. Kazhdan, {\it Arithmetic varieties and their fields of quasi-definition}, Actes du
    Congr$\grave{\rm e}$s International Math\'ematiciens Vol II, Gauthier-Villas, Paris, 1970, 321--325.

\bibitem{Kobayashi83} S. Kobayashi and C. Horst, {\it Topics in complex differential geometry}, DMV Sem. {\bf 3}
    (1983) 4--66.

\bibitem{L-S84} P. Li and R. Schoen, {\it $L^p$ and mean value properties of subharmonic functions on Riemannian manifolds}, Acta. Math. {\bf 153} (1984) 279--301.

\bibitem{Loi06} A. Loi, {\it Calabi's diastasis function for Hermitian symmetric spaces}, Differential Geometry and its Applications {\bf 24} (2006) 311--319.

\bibitem{Mumford79} D. Mumford, {\it An algebraic surface with $K$ ample, $K^2=9$, $p_g=q=0$}, Amer. J. Math.
    {\bf 101} (1979), 233--244.

\bibitem{Ohsawa09} T. Ohsawa, {\it A remark on Kazhdan's theorem on sequences of Bergman metrics}, Kyushu
    J. Math. {\bf 63} (2009), 133--137.

\bibitem{OhsawaTakegoshi87} T. Ohsawa and K. Takegoshi, {\it On the extension of $L^2$holomorphic
    functions}, Math. Z. 195 (1987), 197--204.

\bibitem{Ramadanov67} I. Ramadanov, {\it Sur une propri\'ete de la fonction de Bergman}, C. R. Acad.
    Bulgare Sci. {\bf 20} (1967), 759--762.

\bibitem{Richberg68} R. Richberg, {\it Stetige streng pseudokonvexa Funktionen}, Math. Ann. {\bf 175} (1968), 257--286.

\bibitem{Rhodes93} J. A. Rhodes, {\it Sequences of metrics on compact Riemann surfaces}, Duke Math. J. {\bf
    72} (1993), 725--738.

\bibitem{V-G-P63} \`{E}. B. Vinberg, S. G. Gindikin and I. I. Pjateck\u{i}-\v{S}apiro, {\it Classification and canonical realization of complex homogeneous bounded domains}, Trudy Moskov. Mat. Ob\v{s}\v{c}. {\bf 12}
(1963), 359--388; Trans. Moscow Math. Soc. {\bf 12} (1963), 404--437.

\bibitem{Xu} X. Wang, {\it Effective very ampleness of the canonical line bundles on ball quotients}, to apear in Journal of Geometric Analysis.

\bibitem{Yau77} S. T. Yau, {\it Calabi's conjecture and some new results in algebraic geometry}, Proc. Nat. Acad. Sci. U.S.A.. {\bf 74} (1977), 1798--1799.

\bibitem{Yeung94} S. K. Yeung, {\it Betti numbers on a tower of coverings}, Duke Math. J. {\bf 73} (1994),
    201-226.

\bibitem{Yeung00a} S. K. Yeung, {\it Very ampleness of line bundles and canonical embeddings of coverings
    of manifolds}, Compositio Math. {\bf 123} (2000), 209--223.

\bibitem{Yeung00b} S. K. Yeung, {\it Effective estimates on the very ampleness of the canonical line bundle
    of locally Hermitian symmetric spaces}, Trans. Amer. Math. Soc. {\bf 353} (2000), 1387--1401.

\end{thebibliography}
\end{document}